\documentclass[11pt,twoside]{article}
\usepackage{fancyhdr}
\usepackage[colorlinks,citecolor=blue,urlcolor=blue,linkcolor=blue,bookmarks=false,breaklinks]{hyperref}
\usepackage{amsfonts,epsfig,graphicx}
\usepackage{afterpage}
\usepackage{amsmath,amssymb,amsthm} 
\usepackage{fullpage}
\usepackage[T1]{fontenc} 
\usepackage{epsf} 
\usepackage{graphics} 
\usepackage{amsfonts,amsmath}
\usepackage[sort,authoryear]{natbib} 
\usepackage{psfrag,xspace}
\usepackage{color,etoolbox}
\usepackage{enumitem}

\usepackage{tikz}

\usepackage{setspace}

\usepackage{mathtools}

\newcommand{\inprod}[2]{\ensuremath{\langle #1 , \, #2 \rangle}}

\newtheorem{theorem}{Theorem} 
 
\newtheorem{lemma}{Lemma}

\newtheorem{assumption}{Assumption}

\newtheorem{example}{Example}

\newcommand{\abscontmom}{\mathcal{P}_{2,\mathrm{ac}}(\Omega)}

\newcommand{\semidualpq}{\mathcal{S}_{P,Q}}

\newcommand{\probspace}{\mathcal{P}_2(\Omega)}
\newcommand{\probspacesquared}{\mathcal{P}_2(\Omega \times \Omega)}

\newcommand{\probone}{P}
\newcommand{\probtwo}{Q}

\usepackage{macros}

\begin{document}

\begin{center} {\LARGE{\bf{
Statistical Inference for Optimal Transport Maps: \\
\vspace{.25cm}
Recent Advances and Perspectives
}}}
\\

\vspace*{.3in}

{\large{
\begin{tabular}{ccccc}
Sivaraman Balakrishnan$^{\dagger}$ & Tudor Manole$^{\ddagger}$ & Larry Wasserman$^{\dagger}$ \\
\end{tabular}

\vspace*{.1in}

\begin{tabular}{ccc}
$^\dagger$Department of Statistics and Data Science \\
Machine Learning Department \\
Carnegie Mellon University \\
\texttt{\{siva,larry\}@stat.cmu.edu} 
\end{tabular}

\vspace*{.2in}

\begin{tabular}{ccc}
$^\ddagger$Statistics and Data Science Center \\
Massachusetts Institute of Technology \\
\texttt{tmanole@mit.edu} 
\end{tabular}
\vspace*{.2in}}}

\vspace*{.2in}

\begin{abstract}
\noindent 
In many applications of optimal transport (OT), the object of primary interest is the optimal transport \emph{map}. This map rearranges mass from one probability distribution to another in the most efficient way possible by minimizing a specified cost. In this paper we review recent advances in estimating and developing limit theorems for the OT map, using samples from the underlying distributions. We also  review parallel lines of work that establish similar results for special cases and variants of the basic OT setup. We conclude with a discussion of key directions for future research with the goal of providing practitioners with reliable inferential tools.
\end{abstract}
\end{center}

\section{Introduction}
Optimal transport (OT) is concerned with the following question: Given two probability measures $\probone$ and $\probtwo$ supported in $\mathbb{R}^d$, how can we transport $\probone$ to $\probtwo$ while minimizing a specified transportation cost? 
One way of answering this question is through the notion of an {\it optimal transport map}. 
A (quadratic) optimal transport map $T_0$ from $\probone$ to $\probtwo$
is any solution to the so-called Monge problem~\citep{monge1781},
\begin{equation}
\label{eqn:monge}
\argmin_{T:T_\# P = Q} \int \|T(x)-x\|_2^2 dP(x),
\end{equation}
where the minimizer is chosen from the set of all
transport maps between $P$ and $Q$, that is, the set of Borel-measurable functions 
$T:\bbR^d \to \bbR^d$ such that $T_{\#} P := P(T^{-1}(\cdot)) = Q$. Intuitively, a transport map $T$ is a transformation
of a random variable~$X$ with distribution~$P$ to a random variable~$T(X)$ with distribution~$Q$.
Among the potentially large collection of such transport maps,
the OT map is the one which is most parsimonious, in the sense that it is the
one which is closest to the identity map, on average. 

Studying the existence, uniqueness and regularity properties of the optimal transport map~$T_0$, and developing numerical methods for computing it, have been the focus of extensive research \citep{villani2003,villani2008,peyre2019computational}. In this review, we primarily focus on the problems of estimating $T_0$ and developing limit laws for it when the underlying distributions $P$ and $Q$ are unknown, based on random samples $X_1, \ldots, X_n \sim P$ and $Y_1, \ldots, Y_m \sim Q$. This question is motivated by a growing
number of   statistical applications of OT maps, some of which we survey in Section~\ref{sec:app}. 
More broadly, performing inference for OT maps is a problem which falls within the scope of \emph{statistical}
optimal transport \citep{chewi2024}.

Adopting a statistical perspective on optimal transport problems has several benefits. 
We formally introduce the minimax framework in Section~\ref{sec:background}. 
The minimax framework provides a foundation for comparing estimators of the OT 
map based on their sample efficiency. This framework also allows the statistician to 
characterize the fundamental limits of estimation from random samples and consequently 
to reason about optimal estimators. Estimators based on random samples are inherently 
stochastic and limit laws describe their fluctuations. Limit laws enable rigorous uncertainty 
quantification -- which is often useful when applying the tools of OT to analyze real data sets. 
The statistical perspective complements mathematical and computational approaches, highlighting distinct yet interconnected concerns. These other perspectives are the subject of many  books~\citep{villani2003,villani2008,rachev1998,santambrogio2015,peyre2019computational,figalli2021invitation}.

Many statistical applications of optimal transport stem from 
an alternate characterization of OT maps, which links them to a notion of {\it multivariate monotonicity}. 
This characterization 
can be understood through a
celebrated result of~\cite{brenier1991}, which we state here 
informally.

\begin{theorem}[Brenier's Polar Factorization Theorem (Informal)]
\label{thm:polar_factorization}
For any absolutely continuous probability measures $P$ and $Q$ on $\bbR^d$, 
there exists a unique
optimal transport map $T_0$ from $P$ to $Q$ which can be written as the gradient of a convex function $\varphi_0$, 
i.e. $T_0 = \nabla \varphi_0$. 
Furthermore, for any other  transport map $T$ from $P$ to $Q$, 
the following decomposition holds
$$T = T_0 \circ S,$$
where $S$ is a   transport map from $P$ onto itself. 
\end{theorem}
\noindent By analogy to univariate convex functions, 
whose derivatives are monotonic,
it is natural that multivariate convex functions
have gradients which satisfy a multivariate notion of monotonicity.  In this language, 
Brenier's polar factorization theorem implies that 
the OT map  
is the {\it unique multivariate monotonic map} between two distributions $P$ and~$Q$.
Furthermore, the theorem shows that any other transport map $T$ from $P$ to $Q$ can be written 
as a composition of this monotonic vector field $T_0$, and a  map $S$  which is 
redundant as far as the problem of transporting $P$ to $Q$ is concerned. This gives an alternative sense in which
$T_0$ is most parsimonious among all transport maps, since it is the one for which 
the redundant component is $S = \mathrm{Id}$. 
These various perspectives suggest that $T_0$ can be viewed as a canonical transformation between multivariate random 
variables. In the special case where the distributions are univariate, the unique monotonic 
map between them is the classical quantile-quantile (QQ) function, 
and OT gives a natural generalization of this object to the multivariate setting.
 
Inspired by the OT framework, and aiming to mitigate some of its computational and statistical deficiencies, many variants have been proposed in the recent years. These variants include, for instance, entropic OT \citep{cuturi2013}, Gaussian-smoothed OT~\citep{goldfeld2020convergence}, sliced OT~\citep{bonneel2015sliced}, linear OT~\citep{flamary2019concentration}, linearized OT~\citep{wang2013linear},
and rectified flow~\citep{liu2023flow}. In this review, we also briefly discuss some recent results studying estimation and limit laws for transport maps in these OT-inspired settings.

The rest of the paper is organized as follows. In Section~\ref{sec:app} we discuss applications of OT in which the transport map is a central object of interest. In Section~\ref{sec:background} we provide some background on optimal transport, and on transport maps. 
In Section~\ref{sec:main} we consider the estimation of and limit laws for multivariate OT maps, reviewing a few of the key ideas from some lines of recent work. In Section~\ref{sec:variants} we discuss variants and special cases of OT, surveying some recent results on estimation and limit laws in these settings. In Section~\ref{sec:discussion} we conclude with a discussion of numerous interesting avenues for future work.

\section{OT Maps as the Target of Inference: Motivating Examples}
\label{sec:app}

Before providing further technical background on the optimal transport problem, we provide a glimpse into the wide variety of recent statistical applications
that rely on OT maps.
Our emphasis is on  applications
where OT maps are the main object of interest and for which inferential methods are an important tool;
further applications are discussed
in the surveys of~\cite{panaretos2019a,panaretos2020,chewi2024,kolouri2017}.

\subsection{Statistical Inference with Unpaired Samples}
 \label{sec:app_uncoupled}
One of the early applications of transport-based methods has its roots in the work of
\citet{degroot1980estimation}. \citet{degroot1980estimation} were broadly interested in statistical estimation
questions given a broken sample, where $n$ paired samples $(X_i,Y_i)$ 
are drawn from a joint distribution $\gamma$ with marginals $P$ and $Q$, but in the observed sample,
the correspondence between pairs is ``broken'' by an unknown permutation. That is, one merely observes
the pairs $(X_i,Y_{\tau(i)})$, for some unknown permutation $\tau$ over $\{1,\dots,n\}$. One canonical estimation task in this setting is to estimate the transformation $\tau$ from the broken sample.
A variant of this setting, known as the unlinked setting~\citep{balabdaoui2021unlinked,slawski2024permuted}, is where the observed samples $X_i$ and $Y_i$ are each drawn independently $P$ and $Q$
respectively,
and the goal is to infer some properties of the joint distribution $\gamma$. Although
$\gamma$  is generally unidentifiable, it can be identified in many useful special cases \citep{pananjady2018linear,carpentier2016learning,collier2016,rigollet2019uncoupled,slawski2024permuted,ghodrati2022distribution}. 
One such special case is when 
the joint distribution $\gamma$ is induced by the OT map $T_0$ from $P$ to $Q$, in the sense that $(X,Y) \sim \gamma$ if and only if $Y = T_0(X)$.  This hypothesis  can be justified in some cases by physical constraints, i.e. roughly that samples from the target distribution are obtained by evolving samples from the source distribution for a small amount of time by an ``energy-minimizing'' vector field~\citep{schiebinger2019,bunne2024optimal}.

We illustrate this setting
through an application to single-cell genomics, 
which was studied by~\cite{schiebinger2019}.
The goal of their work was to longitudinally model the  evolution of cell populations 
in biological tissues, 
using single-cell sequencing technologies. 
In a simplified version of their
setup, one assumes
that at any time
$t$, a  population of cells (modeled as $d$-dimensional feature
vectors)
consists of independent 
draws from a probability
distribution $P_t$ on $\bbR^d$. 
The joint distribution of any given
cell at two times $t < s$ is denoted by~$\gamma_{ts}$, with marginal
distributions $P_t$ and $P_s$. 
The 
goal is to estimate various
properties of these joint distributions.
Due to the destructive nature
of single-cell sequencing
technologies, one cannot
sequence a cell 
at more than one timepoint, and
thus one does not have access
to paired observations drawn from
the joint distributions $\gamma_{ts}$.
Instead, we have access to 
samples from $P_t$
at multiple timepoints $t$,
which are {\it independent} across time.
Although the joint
distributions $\gamma_{ts}$ 
are unidentifiable under
this sampling model, 
\cite{schiebinger2019}
postulate the so-called {\it optimal transport principle}, namely
that for all sufficiently
close timepoints $t  < s$, 
$$(X^t,X^s) \sim \gamma_{ts}~~\Longrightarrow ~~ X^s \approx T_{ts}(X^t),$$
where $T_{ts}$ is the OT map
from $P_t$ to $P_s$. 
Under this modeling assumption, 
the temporal evolution
of cells is entirely governed
by the OT maps $T_{ts}$, 
and an important statistical question
is to perform inference
for these maps on the basis of
i.i.d. samples at each timepoint.

\subsection{Measures of Spatial Association}
\label{sec:app_OTC}
 
Optimal transport has also been used to define
new measures of association or dependence between multivariate distributions. 
On the one hand, optimal transport provides a meaningful approach for
extending classical univariate rank-based coefficients of association to the multivariate 
setting~\citep{chatterjee2022}. 
The key to this extension is the fact that optimal transport yields a notion 
of {\it multivariate ranks}, as we discuss further in Section~\ref{sec:app_methods} below. 
On the other hand, a variety of optimal transport divergence functionals
can be used to construct measures of dependence which are similar in spirit to the mutual 
information~\citep{xiao2019,mori2020,mordant2022,nies2021,wiesel2022}. 
In what follows, we specifically   highlight how optimal transport maps
have   been used to 
provide new measures of {\it spatial} association in applications
to the physical sciences~\citep{tameling2021,komiske2019,naas2024multimatch}.  

A  classical measure of association is Ripley's~$K$-statistic \citep{ripley1976}. 
In the two-sample setting, Ripley's $K$-statistic is an unbiased estimator of the functional
$$\calK_\tau (P,Q) = \bbP\big(\|X-Y\|_2 \leq \tau\big),\quad \tau \geq 0,$$
where $X$ and $Y$ are independent random variables, drawn
respectively from distributions $P$ and $Q$ on $\bbR^d$. 
The quantity $\calK_\tau (P,Q)$ is a measure of spatial homogeneity between the two distributions, 
which can be used to gauge their typical level of overlap at varying scales $\tau$. 

Ripley's coefficient is used in the field of biophysics,
to quantify the spatial proximity (also known as {\it colocalization}) 
between groups of biological molecules (e.g.~\cite{lagache2013}). 
Quantitative measures of this type provide
a proxy for the level of interaction between different
types of molecules, and 
indicate whether they are jointly involved in any relevant biological
pathways.
In this setting, the measures $P$ and $Q$ typically represent the densities of two different types
of molecules over a two-
or three-dimensional space. Ripley's $K$-statistic
can either be used as a summary statistic of the colocalization of the molecules, or as
a test statistic for determining whether there is significant interaction between the two groups. 
In either case, $\calK_\tau(P,Q)$ is a quadratic functional which can be estimated with a U-statistic, 
and statistical inference for the $K$-statistic is a routine task~\citep{hoeffding1948class,serfling1980approximation}.

\cite{tameling2021} recently proposed an alternative coefficient, 
called the {\it Optimal Transport Colocalization} (OTC) curve, which is defined by:
\begin{equation}
\label{eq:OTC}
\calO_\tau(P,Q) = \bbP\big( \|X-T_0(X)\|_2 \leq \tau\big), \quad \tau \geq 0,
\end{equation}
where $T_0$ is the optimal transport map pushing $P$ forward onto $Q$.
Much like Ripley's coefficient, the OTC curve aims to compare the typical distance between $P$ and $Q$
over a continuum of scales $\tau$, however the differences being compared are fundamentally different. 
Whereas Ripley's coefficient measures the typical distance between independent random variables $X$ and~$Y$, 
the OTC coefficient $\mathcal{O}_\tau$ chooses $X$ and $Y=T_0(X)$ to be maximally dependent. In the context of colocalization, one can think of $T_0$ 
as  representing the interactions between one group
of molecules with the other, and the OTC curve measures the typical
distance at which such molecules interact. While the OT map is not the only transport map
that could be used in this definition, it is a natural default choice 
if the underlying objects are expected to interact according to a minimum-action principle. 
Another important feature of the OTC coefficient is that it is robust to convolutional noise in the underlying distributions, which arises due to the limited resolution of the imaging devices used to observe the molecules~\citep{tameling2021}. Indeed, many practical applications of OT benefit from its robustness to small amounts of convolutional noise.

It is natural to 
ask whether there exist inferential methods for the OTC coefficient when
the practitioner is merely given samples from $P$ and $Q$. 
This question arises in the biophysical applications described above, 
since many  microscopy technologies are photon-limited, 
thus the images they produce are subject to stochastic photon-counting noise~\citep{aspelmeier2015modern}. 
Since the OTC curve depends on evaluations of the unknown optimal transport map $T_0$ over its
entire support, one is led to the question of performing inference for the OT map itself.

Beyond the OTC a variety of other OT-based  functionals can be used to measure spatial association,
owing to their ability to accurately capture geometric features of the underlying distributions.
Another example
arises in 
high energy physics,
where OT is used to quantify the spatial proximity of energy deposits produced by different  collisions
within particle colliders~\citep{komiske2019,manole2024background,ba2023}. 
When such measures must be estimated from data, inferential methods for OT maps and their functionals are useful tools.

\subsection{Discretization of Numerical Procedures}
 \label{sec:app_discretization}

A third setting which motivates 
  inferential 
methods for OT maps
is the problem
of quantifying the error
induced by numerical methods
for solving the OT problem~\citep{sommerfeld2019optimal}. 
To elaborate, 
the Monge optimal transport 
problem~\eqref{eqn:monge}
is an infinite-dimensional
optimization problem
which is rarely solvable
in closed-form, and essentially
all existing numerical
methods for 
approximating the optimal map $T_0$
involve some form of discretization
of the absolutely continuous measures
$P$ and $Q$. 
For example, one
class of methods
are semi-discrete solvers~\citep{merigot2011,kitagawa2019}, 
which are based on discretizing
the target measure $Q$; 
such solvers are used
in a variety of applications 
such as in  meteorology~\citep{bourne2022semi}
and 
cosmology~\citep{levy2021fast}.
Other examples include the gradient-descent based
methods of~\cite{chartrand2009gradient,jacobs2020fast}, which use
evaluations of the densities of $P$ and $Q$
over a grid.  

Considering the semi-discrete setting
as an example, the
aforementioned numerical methods 
provide efficient approximations of the OT map
$T_n$ between $P$ and any discrete
measure $Q_n$ supported
on $n$ points. The measure $Q_n$ is a user-specified approximation of the absolutely
continuous measure $Q$. 
A natural question is whether
one can quantify the gap
between the OT map $T_n$ from $P$
to $Q_n$ (which is approximated 
by the algorithm), and the true OT map
$T_0$ from $P$ to $Q$. 
When the discretized measure $Q_n$
is an empirical measure comprised
of points sampled from the measure $Q$, 
the inferential methods
surveyed in this paper provide
upper bounds on 
the approximation error of $T_n$ as an
estimator of $T_0$, and in some
cases provide limit laws which sharply
characterize this error in various senses.

\subsection{Monge-Kantorovich Ranks and Quantiles} 
 \label{sec:app_methods}

Another application of optimal transport stems from
the work of~\cite{chernozhukov2017,hallin2021}, who  
observed that optimal transport can be used to 
generalize traditional 
univariate statistical concepts, such as quantile functions, cumulative
distributions functions (CDFs), and their
empirical analogues, based on ranks and signs.   
In order to make this connection clear, let us
recall the classical probability integral transform, 
which states that the  quantile function $G^{-1}$ of  a probability
distribution $Q$ on $\bbR$ is a transport map from the uniform 
distribution $P = \calU[0,1]$ onto $Q$. 
Under mild conditions, the map $G^{-1}$ is the OT map which pushes $P$ onto $Q$. 
This observation suggests a very natural way of
defining {\it multivariate} quantiles: Given a measure $Q$ supported
in $\bbR^d$, and given a reference measure $P$, such as the uniform distribution
on $[0,1]^d$, one {\it defines} the quantile function of $Q$ as the unique optimal transport map $T_0$ pushing $P$ forward onto $Q$. This map is referred
to as the {\it Monge-Kantorovich quantile function} of $Q$. 
Notions of cumulative distribution
functions, ranks, and signs, can also be defined using a similar principle. 

\cite{chernozhukov2017} and \cite{hallin2021} 
argued that these definitions
retain many desirable properties of their univariate  
counterparts. 
Follow-up work has then used these definitions to propose multivariate rank-based, distribution-free 
procedures for testing hypotheses of goodness-of-fit, two-sample equality, independence, 
and symmetry, which often retain the same asymptotic efficiency properties 
as their univariate counterparts (cf.~\cite{shi2020,deb2021a, deb2021b,huang2023}). 
We refer to~\cite{chatterjee2022, hallin2022} for  surveys of these developments. 
Constructing estimators and confidence intervals or bands for
CDFs and quantiles is one of the most
basic tasks in classical nonparametric statistics, and adapting
these ideas to Monge-Kantorovich CDFs and quantiles
requires reasoning about statistical inference for OT maps~\citep{ghosal2022}. 
As a concrete example, constructing a confidence interval for the multivariate median of a distribution requires constructing a pointwise confidence interval for the map $T_0$, 
a topic we return to in Section~\ref{sec:clt_smooth}.

\section{Background and Problem Setting}
\label{sec:background}
The OT problem in~\eqref{eqn:monge} does not have a solution in general, and the pushforward constraint is typically non-linear. 
This motivates the study of the Kantorovich relaxation \citep{kantorovich1948,kantorovich1942}
of identifying an \emph{optimal coupling} between $P$ and $Q$, defined~as: 
\begin{align}
\label{eqn:kantorovich}
\pi_0 \in \argmin_{\pi \in \Pi(P,Q)} \int \|x-y\|_2^2 d\pi(x,y).
\end{align}
The set $\Pi(P,Q)$ denotes the set of couplings of $P$ and $Q$, i.e.
\begin{align*}
\Pi(P,Q) = \{ \pi \in \probspacesquared: \pi( \cdot \times \Omega) = P, \pi( \Omega \times \cdot) = Q \}.
\end{align*} 
Here we let $\Omega \subseteq \mathbb{R}^d$ be a convex set which denotes the domain of the distributions $P$ and $Q$, and we use $\calP(\Omega)$ to denote the set
of Borel probability measures supported
on $\Omega$, and $\probspace$ be the subset of such probability measures with finite second moment. The optimal \emph{value} of the Kantorovich program with quadratic cost defines a metric on $\probspace$ which is known as the $2$-Wasserstein distance:
\begin{align}
\label{eq:w2}
    W_2^2(P,Q) := \min_{\pi \in \Pi(P,Q)} \int \|x - y\|_2^2 d\pi(x,y).
\end{align}
In contrast to the Monge problem~\eqref{eqn:monge}, the Kantorovich problem~\eqref{eqn:kantorovich} is always feasible, and is a linear program. The Kantorovich problem is a relaxation of the Monge problem in the sense that every valid transport map $T$ between $P$ and $Q$ defines a coupling $\pi = (\mathrm{Id},T)_\# P$ between them which achieves the same cost for both objectives.

\subsection{Duality and Correspondences}
The Kantorovich problem in~\eqref{eqn:kantorovich} being a linear program has a dual linear program. Indeed, Kantorovich introduced linear programming duality primarily as a tool to study his OT program \citep{kantorovich1942}. 
In our review, the so-called semi-dual problem  which is derived from the Kantorovich dual will be of central interest. The Wasserstein distance has a characterization via the identity:
\begin{align*}
W_2^2(\probone, \probtwo) = \int \| \cdot\|_2^2 dP + \int \| \cdot\|_2^2 dQ - 2 \inf_{\varphi \in L^1(P)} \semidualpq(\varphi),
\end{align*}
where the semi-dual functional $\semidualpq$ is:
\begin{align*}
    \mathcal{S}_{P,Q}(\varphi) = \int \varphi dP + \int \varphi^* dQ,
\end{align*}
and $\varphi^*(y) := \sup_{x \in \Omega} \left[\inprod{x}{y} - \varphi(x)\right]$,
denotes the Fenchel conjugate of $\varphi.$
The semi-dual problem is to optimize the semi-dual functional:
\begin{align}
\label{eqn:semidual}
     \varphi_0 := \argmin_{\varphi \in L^1(P)} \mathcal{S}_{P,Q}(\varphi).
\end{align}
One link between the Monge and Kantorovich problems is given by Brenier's theorem, a version
of which was already stated in the Introduction. Brenier's theorem shows that, under the assumption that the source measure is absolutely continuous, i.e. $P \in \abscontmom$, the solutions to these problems are isomorphic:
\begin{theorem}[Brenier's Theorem]
\label{thm:brenier}
Let $P \in \abscontmom$ and $Q \in \probspace$. 
Then, the following assertions hold.
\begin{enumerate}
\item  
There exists an optimal transport map $T_0$ 
pushing $P$ forward onto $Q$
which takes the form $T_0 = \nabla\varphi_0$ for a convex function $\varphi_0: \mathbb{R}^d \to \mathbb{R}$
which solves the semi-dual problem~\eqref{eqn:semidual}. Furthermore, $T_0$
  is uniquely determined $P$-almost everywhere.  
\item If we further have $Q \in \abscontmom$, then $S_0 := \nabla \varphi_0^*$ is 
the $Q$-almost everywhere uniquely
determined optimal transport map pushing $Q$ forward onto $P$. 
\end{enumerate}
\end{theorem}
\noindent When Brenier's theorem holds there is a direct relationship between the solutions to Monge's problem and Kantorovich's problem. An important consequence of this tight connection is that we can use tools from convex analysis and convex duality, which apply directly to Kantorovich's problem, to study estimation and limit laws for OT maps. From a methodological standpoint, as we discuss further in Section~\ref{sec:estimators}, the dual perspective will also suggest a class of natural estimators of the transport map which are obtained from approximate solutions to the semi-dual problem~\eqref{eqn:semidual}.

\subsection{Special Cases}
The optimal transport problem greatly simplifies in certain special cases. When the distributions $P$ and $Q$ are univariate, then the OT map has a simple closed-form expression. Suppose that $P$ has no atoms. Then, the OT map from $P$ to $Q$ is given by
the classical quantile-quantile function
\begin{align*}
    T_0(x) = F_{Q}^{-1}(F_P(x)),
\end{align*}
where $F_P(x) := P(X \leq x)$ is the cumulative distribution function of $P$, and $F_Q^{-1}$ is the quantile function of $Q$ (see, for instance, \cite[Chapter 2]{santambrogio2015}). 
In this case, a precise characterization of various statistical OT quantities can be obtained 
\citep{bobkov2019}, based on
classical results for quantile and empirical processes~\citep{shorack1986empirical}. This one-dimensional theory can be lifted to multivariate settings, while retaining some of the one-dimensional benefits in various ways~\citep{bonneel2015sliced, rosenblatt1952remarks,knothe1957contributions,paty2019subspace}.

If $P = N(\mu_1, \Sigma_1)$ and $Q = N(\mu_2,\Sigma_2)$ are multivariate Normal, then 
once again the OT map from $P$ to $Q$ has a closed form expression:
\begin{align}
\label{eqn:gaussianot}
    T_0(x) = \Sigma_1^{-1/2} (\Sigma_1^{1/2} \Sigma_2 \Sigma_1^{1/2})^{1/2} \Sigma_1^{-1/2} (x -\mu_1) + \mu_2.
\end{align}
The transport map in this setting depends only on the means and covariances of $P$ and $Q$. As a consequence, to study estimation and limit laws for the transport map, we can directly build on the extensive literature on multivariate mean and covariance estimation. 
More generally, suppose we consider general distributions $P$ and $Q$ which have finite first and second moments (not necessarily multivariate Gaussian), but require that $Q = T_{\#} P$, where $T$ is affine. Then the OT map from $P$ to $Q$ once again has the same closed form expression~\eqref{eqn:gaussianot} \citep[Lemma 1]{flamary2019concentration}.

\subsection{Entropic Optimal Transport}
The entropic optimal transport problem corresponds to a regularized version of Kantorovich's problem~\eqref{eqn:kantorovich} and has its roots in the work of \citet{Schrodinger1931} (see \citep{leonard2014survey}). In entropic OT, the goal is to identify a coupling which solves the following program:
\begin{align}
\label{eqn:eot}
\pi_{\varepsilon} \in \argmin_{\pi \in \Pi(P,Q)} \int \|x - y\|_2^2 d\pi(x,y) + \varepsilon \text{KL}(\pi \| P \otimes Q),
\end{align}
where the Kullback-Leibler regularizer encourages the coupling $\pi_{\varepsilon}$ to be ``spread out'' relative to the OT coupling $\pi_0$. 

The computational benefits of the entropic formulation were highlighted in the work of \citet{cuturi2013}, who showed that this problem could be solved (in the discrete case) to reasonable accuracy even for large-scale problems via the Sinkhorn algorithm 
(also known as iterative proportional fitting).
This work set the stage for a resurgence of interest in numerical methods for solving OT problems, and many of its variants~\citep{peyre2019computational}.
The convergence of $\pi_{\varepsilon}$ to $\pi_0$ as $\varepsilon \rightarrow 0$ is also well-studied (see, for instance, \citep{nutz2022entropic,carlier2023convergence}). In fact, as $\varepsilon \rightarrow 0$, the induced measure of divergence converges to the Wasserstein distance, and various dual objects for the two optimization problems~\eqref{eqn:eot} and~\eqref{eqn:kantorovich} converge as well, making entropic OT a useful tool in the study of its classical, unregularized counterpart. In addition, the entropic OT map with fixed $\varepsilon > 0$ also has its own distinct interpretation as an ``energy-minimizing'' coupling, and in some scientific applications the entropic OT map is a natural target of inference in its own right (see~\citep{rigollet2025sample}).

The entropic OT problem has a dual program which can be viewed as a relaxation of Kantorovich's dual for the OT problem:
\begin{align}
\label{eqn:eotdual}
    (f_\varepsilon, g_\varepsilon) := \argmax_{f,g \in \mathcal{C}_b(\Omega)} \int f dP + \int g dQ - \varepsilon \iint \left[\exp((f(x) + g(y) - \|x - y\|_2^2)/\varepsilon) - 1\right] dP(x) dQ(y),
\end{align}
where $\mathcal{C}_b(\Omega)$ is the set of bounded continuous functions supported on $\Omega.$ 
For the entropic OT problem, one could define the ``entropic map'' as \citep{pooladian2021}:
\begin{align*}
    T_{\varepsilon}(x) := \mathbb{E}_{\pi_{\varepsilon}}[Y | X = x],
\end{align*}
which is a projection of the entropic coupling $\pi_{\varepsilon}$ onto the space of maps. The entropic map $T_{\varepsilon}$ as well as entropic coupling $\pi_{\varepsilon}$ have simple expressions in terms of the dual potentials $f_{\varepsilon}$ and $g_{\varepsilon}$ \citep{pooladian2021}. It is important to note that in general $T_{\varepsilon}$ is not a valid transport map between $P$ and $Q$, i.e. $T_{\varepsilon\#} P \neq Q$.

When the distributions under consideration are supported on a compact set $\Omega$, the entropic OT problem satisfies some very strong regularity properties. For instance, the solutions to the dual program to~\eqref{eqn:eot} have all derivatives bounded (see, for instance, \citep[Lemma 1]{goldfeld2024limit}). These strong regularity properties underlie the many statistical benefits of entropic OT, wherein estimation and limit laws for the entropic map can be obtained at parametric rates \citep{rigollet2025sample,goldfeld2024limit,gonzalez-sanz2022weak,genevay2019sample,mena2019statistical,gunsilius2021matching}. 

\subsection{Statistical Setup}
In the statistical OT setting, we observe samples $X_1,\ldots,X_n \sim P$ and $Y_1,\ldots,Y_m \sim Q$, and our goal is to estimate the OT map $T_0$ between $P$ and $Q$. Following~\citet{hutter2021} one way to evaluate an estimate $\widehat{T}$ is by its risk:
\begin{align}
\label{eqn:risk}
\bbE\|\hat T - T_0\|_{L^2(P)}^2 = \mathbb{E} \left[\int \|\widehat{T}(x) - T_0(x)\|_2^2 dP(x) \right],
\end{align}
where the outer expectation is taken over the randomness in the samples. Letting $\widetilde{Q} := \widehat{T}_{\#}P$, we observe that $L^2(P)$-loss function used above measures a distance between $\widetilde{Q}$ and $Q$ which upper bounds the usual Wasserstein distance, i.e.:
\begin{align*}
    W_2^2(\widetilde{Q},Q) \leq \int \|\widehat{T}(x) - T_0(x)\|_2^2 dP(x).
\end{align*}
The quantity on the right-hand side above is sometimes
known as the linearized Wasserstein distance~\citep{wang2013linear} where $P$ is used as a reference measure. The linearized Wasserstein distance is often used as an easy-to-compute proxy for the Wasserstein distance (e.g.\,\cite{cai2020linearized}). 

When it comes to limit laws for the transport map, and more broadly to uncertainty quantification, several quantities are of interest in applications. For instance, one goal is to study the pointwise fluctuations of the transport map. Here for a fixed $x_0 \in \Omega$ the goal is to characterize the limiting distribution $W$ of the appropriately rescaled error \citep{manole2023}:
\begin{align*}
r_n (\widehat{T}(x_0) - T_0(x_0)) \cdist W,
\end{align*}
where $r_n$ quantifies the rate of convergence of $\widehat{T}(x_0)$ to $T_0(x_0)$.
For constructing confidence bands it is useful to instead study the distribution of errors uniformly over the domain~\citep{ponnoprat2024uniform}:
\begin{align*}
r_n \sup_{x \in \Omega} \frac{\widehat{T}(x) - T_0(x)}{\sigma_x} \cdist W,
\end{align*}
for a scale parameter $\sigma_x$. 
Weak limits, which quantify the limiting distribution of the entire process $r_n (\widehat{T} - T_0)$ 
have also been studied (particularly for the entropic OT map) \citep{gonzalez-sanz2022weak,goldfeld2024limit,gonzalez2023weak,harchaoui2020asymptotics,gunsilius2021matching}. 
These weak limits allow us to quantify the uncertainty of the transport map or coupling as measured via a 
``test function'' $\eta$. For instance, denoting by $\widehat{\pi}_{\varepsilon}$ an estimate 
of the entropic OT coupling, these papers quantify the limiting distribution $W$ of:
\begin{align}
\label{eq:test_function_clt}
    r_n \int \eta(x,y) (\widehat{\pi}_{\varepsilon}(x,y) - \pi_{\varepsilon}(x,y)) \cdist W.
\end{align}
As an illustration, suppose we fix two sets $A, B \subseteq \Omega$, and our goal is to produce an asymptotically valid confidence set for the fraction of mass from the set $A$ which is coupled to points in the set $B$ -- then we could apply the above result with the test function $\eta(x,y) = \mathbb{I}(x \in A, y \in B)$. A special case of this functional is precisely the optimal transport colocalization
curve $\calO_\tau$ at a fixed location $\tau$, described in Section~\ref{sec:app_OTC}. As a second example, 
one can apply the above result to the test function $\eta(x,y) = \|x-y\|^2$ to obtain the (regularized) Wasserstein 
distance\footnote{Distributional limits and minimax estimation rates for this particular 
functional have been developed in some generality;  see for instance~\cite{del2023,hundrieser2024unifying,hundrieser2024lca,manole2024sharp}
and references therein.}.

Beyond establishing the existence of limit laws, establishing the validity
of the bootstrap or constructing reasonable estimates of the limiting distribution are important steps which enable practical inference.

\subsection{Estimators of the OT and Entropic Maps}
\label{sec:estimators}
There are many ways to construct estimates of the OT and entropic maps from the observed samples, and we discuss a few classes of estimators here.

One natural strategy to construct an estimate of the transport map is to first construct an estimate of the optimal coupling between the empirical measures
\begin{equation} 
\label{eq:empirical_measures}
P_n = \frac 1 n \sum_{i=1}^n \delta_{X_i}, \quad \text{and}\quad Q_m = \frac 1 m \sum_{j=1}^m \delta_{Y_j}.
\end{equation}
This amounts to solving Kantorovich's linear program~\eqref{eqn:kantorovich} with $P$ and $Q$ replaced by $P_n$ and $Q_n$, which can be done by the Hungarian algorithm or the Auction algorithm~\citep{kuhn1955hungarian,bertsekas1979distributed}. 
If $m=n$ then this amounts to finding the permutation $\pi$
that minimizes
$n^{-1}\sum_i \|X_i - Y_{\pi(i)}\|^2$.
This does not however yield an estimate of the transport map defined outside the sample points $X_1,\ldots,X_n$. We can obtain a transport map defined on the entire domain using ideas from nonparametric regression (for instance by nearest neighbors or local smoothing) \citep{manole2024plugin}.

An alternative, closely related, class of estimators can be constructed via the plugin principle. Here we first construct estimates $\widehat{\probone}_n$ and $\widehat{\probtwo}_m$ of the measures $\probone$ and $\probtwo$, and then define $\widehat{T}_{nm}$ to be the optimal transport map between $\widehat{\probone}_n$ and $\widehat{\probtwo}_m$  \citep{manole2024plugin,deb2021,gunsilius2021}. For $\hat T_{nm}$ to be well-defined, the estimators
$\hat P_n$ and $\hat Q_m$ need to be {\it proper}, in the sense that they define probability measures
in their own right, and are related by an optimal transport map 
(which is, for instance, the case if $m=n$ for the
empirical distributions, or if $\widehat P_n$ is absolutely continuous, by Theorem~\ref{thm:brenier}). 
In cases where $\hat P_n$ and $\hat Q_m$ define probability measures but are not related by an OT map, one can instead compute an OT coupling $\widehat{\Pi}_{nm}$ and project this coupling onto the space of maps. This projection, known as the barycentric projection, is simply the conditional expectation of the OT coupling:
\begin{align}
\label{eqn:barycenter}
\widehat{T}_{nm}(x) = \mathbb{E}_{\widehat{\Pi}_{nm}}[Y | X = x].
\end{align} 

A third class of estimators is inspired by the semi-dual optimization problem~\eqref{eqn:semidual}. Under the conditions of Brenier's theorem, the OT map $T_0 = \nabla \varphi_0$, where $\varphi_0$ is defined in~\eqref{eqn:semidual}. Taking inspiration from the vast literature on M-estimation~\citep{vandervaart1996,vandegeer2000} a natural class of estimators for $\varphi_0$, is based on solving an empirical counterpart of~\eqref{eqn:semidual}:
    \begin{align}
    \label{eqn:dual}
        \widehat{\varphi} := \argmin_{\varphi \in \Phi} \int \varphi dP_n + \int \varphi^* dQ_m,
    \end{align}
where $\Phi$ is an appropriate function class.
    Under some regularity conditions, solving the empirical semi-dual optimization problem yields an estimate of the transport map $\widehat{T}_{nm} = \nabla \hat{\varphi}$ \citep{hutter2021,divol2022a,vacher2021convex,ding2024}.

In the entropic setting, one natural estimator comes from the work of~\citet{pooladian2021}.
They propose to solve the empirical dual problem~\eqref{eqn:eotdual} where $P$ and $Q$ are replaced by their empirical counterparts $P_n$ and $Q_m$. This is similar to the dual estimators of the OT map described above, but due to the regularity of entropic OT problem we do not need to restrict the dual optimization to be over a restricted class of functions. The unrestricted dual optimization can also be carried out with the Sinkhorn algorithm~\citep{cuturi2013}. The Sinkhorn algorithm yields a pair of potentials $(\widehat{f}_{\varepsilon},\widehat{g}_{\varepsilon})$ which can be smoothly extended to all of $\mathbb{R}^d.$ Then the estimate of entropic OT map is given by:
\begin{align}
\label{eqn:entropicmap}
    \widehat{T}_{\varepsilon}(x) := \frac{\sum_{i=1}^m Y_i \exp ((\widehat{g}_{\varepsilon}(Y_i) - \|x - Y_i\|_2^2)/\varepsilon)}{\sum_{i=1}^m \exp 
    ((\widehat{g}_{\varepsilon}(Y_i) - \|x - Y_i\|_2^2)/\varepsilon)}.
\end{align}

\section{Estimation and Inference for Optimal Transport Maps}
\label{sec:main}
Recent work has begun to study minimax rates for estimating transport maps in various settings \citep{hutter2021,divol2022a,pooladian2023minimax,manole2024plugin,deb2021}. In this section, we first discuss the basic principles that underlie the analysis of the dual estimator in~\eqref{eqn:dual}, which is based on solving an empirical approximation of the semi-dual optimization problem. We then turn our attention to plugin estimators. 

\subsection{Estimating Maps via the Semi-Dual}
It is clear that under the conditions of Brenier's theorem (Theorem~\ref{thm:brenier}), the optimum of the population semi-dual objective $\varphi_0$ is unique (up to constant shifts) and can be directly transformed into the transport map $T_0 := \nabla \varphi_0$.
However, in order to provide meaningful guarantees for the optimizer of the \emph{empirical} approximation of the semi-dual objective we need to 
quantify the stability of the minimizer to perturbations of the objective function. Toward that goal, it is useful to understand the behavior of the semi-dual objective as we move away from its minimizer. The key result in this context comes from the work of \citet{hutter2021}. We state a version from the recent paper of \citet[Lemma 1]{balakrishnan2025stability}:
\begin{lemma}
    \label{lem:bregman}
    Suppose that $P = (\nabla \varphi_0)_{\#}Q$, and let $\varphi$ be a differentiable, strictly convex, function.
    Then:
    \begin{align}
    \label{eqn:bregman}
        \mathcal{S}_{P,Q}(\varphi) - \mathcal{S}_{P,Q}(\varphi_0) = \int \left[ \varphi^*(T_0(x)) - 
\varphi^*(\nabla \varphi(x)) - \inprod{\underbrace{\nabla \varphi^*(\nabla \varphi(x))}_{=x}}{T_0(x) - \nabla \varphi(x)}\right] d P(x). 
    \end{align}
\end{lemma}
\noindent This lemma shows that the growth of the semi-dual function as we move away from the minimizer $\varphi_0$ to a strictly convex function $\varphi$ is exactly captured by a (integrated) Bregman divergence between $T_0$ and $\nabla \varphi$. If we assume that $\varphi$ is not merely strictly convex, but $\alpha$-strongly convex and $\beta$-smooth (see~Assumptions $\text{A}1(\alpha)$, $\text{A}2(\beta)$ in \citep{balakrishnan2025stability}) then we directly obtain from Lemma~\ref{lem:bregman} that \citep{hutter2021}:
\begin{align}
\label{eqn:quad_growth}
    \frac{1}{\beta} \|\nabla \varphi - T_0\|_{L^2(P)}^2 \leq \mathcal{S}_{P,Q}(\varphi) - \mathcal{S}_{P,Q}(\varphi_0) \leq   \frac{1}{\alpha} \|\nabla \varphi - T_0\|_{L^2(P)}^2.
\end{align}
This fact highlights that the semi-dual functional can be a reasonable proxy for the transport map estimation risk~\eqref{eqn:risk}. If we are able to find a smooth function $\varphi$ such that $\mathcal{S}_{P,Q}(\varphi) - \mathcal{S}_{P,Q}(\varphi_0)$ is small then~\eqref{eqn:quad_growth} guarantees us that $\nabla \varphi$ will be a good estimate of the transport map $T_0.$

To illustrate the main ideas we consider the following simplified setting of ``transport map selection'' \citep{vacher2021convex}, where we are given a finite collection of potentials and are tasked with selecting a good one. We make the following assumptions:
\begin{assumption}
    \begin{enumerate}
        \item $P$ and $Q$ are supported on a compact convex set $\Omega$ and $Q = (\nabla \varphi_0)_{\#}P$.
        \item The potential $\varphi_0 \in \Phi := \{\varphi_1,\ldots,\varphi_M\}$. For every $\varphi \in \Phi$, $x,y \in \Omega$, $|\varphi(x)|, |\varphi^*(y)| \leq B$.
        \item The potentials $\varphi_1,\ldots,\varphi_M$ are $\beta$-smooth on $\Omega$.
    \end{enumerate}
\end{assumption}

\noindent Under these assumptions, it is straightforward to analyze the estimator in~\eqref{eqn:dual}. We sketch the analysis here.
We begin with the so-called basic inequality which is a direct consequence of~\eqref{eqn:dual}:
\begin{align*}
S_{P_n,Q_m}(\widehat{\varphi}) \leq S_{P_n,Q_m}(\varphi_0).
\end{align*}
Adding and subtracting some terms we obtain that,
\begin{align*}
S_{P,Q}(\widehat{\varphi}) - S_{P,Q}(\varphi_0) &\leq S_{P_n,Q_m}(\varphi_0) -  S_{P,Q}(\varphi_0) + S_{P,Q}(\widehat{\varphi}) - S_{P_n,Q_m}(\widehat{\varphi}) \\
&\leq 2 \sup_{j \in \{1,\ldots,M\}} |S_{P,Q}(\varphi_j) - S_{P_n,Q_m}(\varphi_j)|. 
\end{align*}
Now, using the smoothness of the potentials~\eqref{eqn:quad_growth}, and via a direct application of Hoeffding's inequality 
we obtain that with probability 
at least $1 - \delta$,
\begin{align*}
\frac{1}{\beta} \|\nabla \widehat{\varphi} - \nabla \varphi_0\|_{L^2(P)}^2 \leq \sqrt{\frac{B\log(M/\delta)}{2n}} + \sqrt{\frac{B\log(M/\delta)}{2m}}.
\end{align*}
This result shows that under our assumptions the semi-dual criterion yields reasonable estimates of the transport map, incurring only 
a mild (logarithmic) dependence on the size of the set of candidates.
The similarities between this analysis and the analysis of M-estimators~\citep{vandegeer2000}, suggest some natural ways in which one can strengthen this analysis 
-- (1) we can study dual potential selection where the class $\Phi$ is infinite using covering and chaining techniques, (2) we can study variants of the estimator in~\eqref{eqn:dual} where the class $\Phi$ is chosen to balance the estimation error which depends on the complexity of $\Phi$ and the approximation error which depends on how closely some member of the class $\Phi$ approximates the true potential $\varphi_0$ and (3) we can sharpen the rates from the slow $1/\sqrt{n}$ rates to fast $1/n$ rates
by exploiting smoothness and strong-convexity via a technique known as localization (our analysis above only uses smoothness of the potentials). These improvements feature prominently in the works of \cite{hutter2021,divol2022a,ding2024} who provide sharper, and in many cases minimax-optimal, analyses of the dual estimator~\eqref{eqn:dual}. Importantly, the work of \citet{divol2022a} provides meaningful guarantees 
for some practically successful dual OT map estimators which use neural networks to parametrize
the class $\Phi$ \citep{makkuva2020,amos2017input}.

\subsection{Stability Bounds}
Before turning our attention to the analysis of concrete plugin estimators we present one possible high-level strategy for analyzing plugin estimators. The optimal transport program~\eqref{eqn:monge} takes as input two distributions $P$ and $Q$ and produces as output the OT map $T_0$ between them. The plugin estimate $\widehat{T}_{nm}$ in~\eqref{eqn:barycenter} is calculated from the same optimization program from the OT coupling between estimates $\widehat{P}_n$ and $\widehat{Q}_m$. It is natural to expect that if $\widehat{P}_n$ were close to $P$ and $\widehat{Q}_m$ were close to $Q$ then it should be the case that $\widehat{T}_{nm}$ is close to $T_0$. Formally justifying this requires identifying when the optimization program~\eqref{eqn:monge} is \emph{stable} to (some) perturbations of its inputs. 
Stability bounds for the OT program have a long history (see the discussion in Section 3.1 of \citet{balakrishnan2025stability}), and further developing our understanding of stability bounds is an active
research area across many disciplines.
Under the conditions of the previous section -- that $Q = (\nabla \varphi_0)_{\#} P$, for a smooth and strongly convex $\varphi_0$ -- a result of this form appears in the recent work of \citet{balakrishnan2025stability}:
\begin{theorem}
\label{thm:stability}
Suppose that $Q = (\nabla \varphi_0)_{\#} P$, for a $\beta$-smooth and $\alpha$-strongly convex function $\varphi_0$. Then for any $\widehat{P}$ and $\widehat{Q}$ the barycentric projection $\widehat{T}$ computed as in~\eqref{eqn:barycenter} satisfies:
\begin{align*}
    \|\widehat{T} - T_0 \|_{L^2(\widehat{P})}^2 \leq \left((1 + \beta) W_2(\widehat{P},P) + \sqrt{\frac{\beta}{\alpha}} W_2(\widehat{Q},Q)\right)^2.
\end{align*}
\end{theorem}
\noindent This theorem highlights that when $\varphi_0$ is smooth and strongly convex the optimal transport program possesses remarkable stability. 
To construct accurate estimates of the transport map, it is sufficient to construct accurate estimates of the distributions $P$ and $Q$ in the Wasserstein sense. This is the familiar statistical task of distribution/density estimation. One caveat is that Theorem~\ref{thm:stability} upper bounds the $L^2(\widehat{P})$ error of our estimate, but the risk we typically study is the $L^2(P)$ error. In concrete applications, which we take up in the next section, we need to relate these two risks.

\subsection{Transport Map Estimation via Nearest Neighbors}
\label{sec:NN}
 
In this section, we consider perhaps the most natural estimator of the transport map -- the nearest neighbor transport map estimate $\widehat{T}_{\text{NN}}$ which is constructed by extending the empirical OT coupling via one-nearest-neighbor extrapolation. Statistical properties of this estimator have been studied in past work by~\citet{manole2024plugin} and \citet{pooladian2023minimax}. This section is based on the recent paper of~\citet{balakrishnan2025stability}.

Let $X_1,\ldots,X_n \sim P$ and $Y_1,\ldots,Y_n \sim Q$ be i.i.d. samples. We take
the two samples to be of the same size in 
order to convey the main ideas, 
however this restriction is not necessary. 
Under this restriction, the
optimal transport problem from
$P_n$ to $Q_n$ is in direct correspondence
to the following
optimal {\it matching} program:
$$
\tau_n \in \argmin_{\tau}
\sum_{i=1}^n  \|X_i-Y_{\tau(i)}\|^2,
$$
where the minimizer is 
taken over all matchings $\tau$, i.e. permutations on 
$\{1,\dots,n\}$.
Any optimal matching $ \tau_n$
induces an optimal transport map
$T_n$ from $P_n$ to $Q_n$, given by
$T_n(X_i) = Y_{\tau(i)}$.
The map
$T_n$ is a natural in-sample
estimator of the optimal transport map, and its in-sample risk has been studied by~\cite{deb2021,manole2024plugin}. One-sample
variants of this estimator have also been studied by~\cite{chernozhukov2017, hallin2021,ghosal2022}.
Under the regularity
conditions of Theorem~\ref{thm:stability}, a direct consequence of Theorem~\ref{thm:stability} is a
bound on the in-sample risk of $T_n$:
\begin{align*}
\|T_n-T_0\|_{L^2(P_n)}^2 =
\frac{1}{n} \sum_{i=1}^n 
     \|T_{n}(X_i) - T_0(X_i)\|_2^2 
     \lesssim W_2^2(Q_n,Q)  + W_2^2(P_n,P).
\end{align*}
We have thus reduced the problem of bounding the risk of $T_n$ to that of bounding the risk of empirical measures under the Wasserstein distance, which is a classical and widely-studied topic~\citep{fournier2015,lei2020,weed2019,boissard2014a,dudley1969}.
Under mild conditions, one has that $\bbE W_2^2(P_n,P) \lesssim n^{-2/d}$
whenever $d \geq 5$, and
this rate is sharp when $P$ and $Q$
are assumed to be absolutely continuous.
We thus obtain the corresponding
upper bound on the in-sample risk of
estimating $T_0$:
\begin{align}
\label{eq:in_sample_risk}
  \bbE\|T_n-T_0\|_{L^2(P_n)}^2
 \lesssim n^{-2/d},\quad \text{for all } d \geq 5.
\end{align}

The estimator $T_n$
is not defined outside of the support
of $P_n$, and thus does not provide
an estimator of $T_0$ over the whole
support of $P$. We obtain
such an estimator by extrapolating $T_n$
using a one-nearest neighbor lookup. 
Concretely, define the Voronoi partition induced by $X_1,\ldots,X_n$ as 
\begin{align*}
V_j = \big\{x \in \bbR^d: \norm{x-X_j} \leq \norm{x-X_i} , \  \forall i\neq j\big\}, \quad j=1, \dots, n.
\end{align*} 
Then, the one-nearest neighbor estimator 
is defined by
\begin{align*}
\widehat T_{n}^{\mathrm{1NN}}(x) = \sum_{i=1}^n   I(x \in V_i) Y_{\tau_n(i)},\quad x \in \bbR^d.
\end{align*}
Let us now briefly describe how
one can bound the {\it out-of-sample}
risk of $\widehat T_n^{\mathrm{1NN}}$. 
Notice that this estimator is a piecewise constant
extension of $T_n$ over the Voronoi
cells, and recall that $T_0$ is assumed to be Lipschitz.
Leveraging these two facts, 
it turns out that one can reduce
the $L^2(P)$ error to the $L^2(P_n)$
error using only upper bounds on the 
typical nearest-neighbor distance, 
and on the maximum  mass
of each Voronoi cell, denoted by:
 \begin{align*}
     \delta := \mathbb{E}_{X \sim P} \left[\|X - \text{nn}(X) \|_2^2 | X_1,\ldots,X_n\right],
    \quad \text{and,} \quad 
\rho := \max_{1\leq i \leq n} P(V_i),
\end{align*}
where $\mathrm{nn}(x)$ denotes
the nearest neighbor of any point $x \in \bbR^d$ among $\{X_1,\dots,X_n\}$.
The connection between these
various quantities is described in the following
simple Lemma. 
\begin{lemma}
\label{lem:nn}
Under the assumptions of Theorem~\ref{thm:stability},
\begin{align*}
\|\widehat{T}_{n}^{\mathrm{1NN}} - T_0\|^2_{L^2(P)} &\lesssim ( n\rho ) \|T_n - T_0\|_{L^2(P_n)}^2 + \delta.
 \end{align*}
\end{lemma}
\noindent 
The quantities $\delta$ and $\rho$ can be controlled under relatively mild moment assumptions on the distributions $P$ and $Q$. Concretely, as shown by \citet{balakrishnan2025stability}, if more than~4 moments of $P$ and $Q$ are bounded, then 
$\mathbb{E}[\delta] \lesssim n^{-2/d}$ and $\mathbb{E}[\rho] \lesssim (\log n)/n$. Putting these results together 
with the in-sample risk bound~\eqref{eq:in_sample_risk}
and with the appropriate formalism, 
we see that under mild moment conditions on the distributions $P$ and $Q$, if $\varphi_0$ is smooth and strongly convex, the nearest neighbor transport map is consistent and has error which scales as:
\begin{align}
\label{eqn:nnrate}
    \mathbb{E} \|\widehat{T}_{n}^{\mathrm{1NN}} - T_0\|^2_{L^2(P)} \lesssim  n^{-2/d}\log n.
\end{align}
Under the conditions we have placed, 
this turns out to be the minimax optimal
rate of estimating $T_0$~\citep{hutter2021}.

It is worth comparing
this result 
to  nonparametric
regression problems 
based on $k$-nearest neighbor
methods, where
one typically
requires a diverging number of neighbors $k$
in order to achieve consistent estimation~\citep{gyorfi2006}.  
Roughly speaking, the choice 
$k=1$
suffices in our context  because,
unlike the nonparametric regression setting,
the bias of the in-sample estimator $T_n$ is
typically larger than its variance~\citep{manole2024sharp}.
Increasing $k$ would only further increase
this bias. 

It is also worth emphasizing, again, that
in addition to mild moment conditions, 
the bound~\eqref{eqn:nnrate}
relies on the smoothness
and strong convexity  of $\varphi_0$. 
We will discuss this assumption in further detail
in Section~\ref{sec:beyond_smooth} below, 
however we present here a simple special case
where these various
conditions are satisfied. 

\begin{example}[Log-Concave Distributions]
Suppose that $P$ and $Q$ are log-smooth and strongly-log concave distributions, meaning
that they are supported on $\mathbb{R}^d$ and admit Lebesgue densities of the form 
$\exp(-V)$ and $\exp(-W)$, where $V$ and $W$ are twice-differentiable, and satisfy for all $x \in \mathbb{R}^d$: 
 \begin{equation}\label{eq:log_concave} \alpha_V I \preceq \nabla^2 V(x) \preceq \beta_V I,\quad \text{and,}\quad \alpha_W I \preceq \nabla^2 W(x) \preceq \beta_W I.
 \end{equation}
Under these assumptions, Caffarelli's contraction theorem~\citep{caffarelli2000} implies that the Brenier potential $\varphi_0$ is  $\alpha$-strongly convex and $\beta$-smooth where $\alpha = \sqrt{\alpha_V/\beta_W}$ and $\beta = \sqrt{\beta_V/\alpha_W}$. Furthermore, the tails of a log-concave distribution are sub-exponential~\citep{ledoux2005concentration}, and these distributions satisfy the moment conditions above, and consequently the nearest neighbor transport map satisfies the risk bound in~\eqref{eqn:nnrate}.
\end{example} 

\noindent The nearest neighbor map is extraordinarily simple and practical. It is completely free of tuning parameters, and is consistent under relatively mild conditions. On the other hand, the nearest neighbor estimator cannot exploit higher amounts of smoothness in the source and target distributions, and we turn our attention to other (smoother) plugin estimators next.

\subsection{Smooth OT Map Estimation}
\label{sec:smooth}
The stability bound of 
Theorem~\ref{thm:stability} can also be used to characterize 
the risk of smooth plugin estimators of the OT map. 
In this section, we continue to focus on 
the setting where $Q = (\nabla \varphi_0)_{\#} P$ for a smooth and strongly 
convex function $\varphi_0$.
Given i.i.d. samples $X_1,\dots,X_n \sim P$ and $Y_1,\dots,Y_m \sim Q$, we 
let $\hat P_n$ and $\hat Q_m$
denote estimators of $P$ and $Q$. We assume that $\hat P_n$ and $\hat Q_m$ are {\it proper}
estimators, in the sense that they    
almost surely define valid probability measures. 
We let $\widehat{T}_{nm}$ denote the barycentric projection of the OT coupling between $\hat P_n$ and $\hat Q_m$. 
This uniquely defined map~$\hat T_{nm}$ is a natural
plugin estimator of $T_0$, and 
Theorem~\ref{thm:stability} implies the
following bound on its $L^2(\hat P_n)$ risk:
$$\bbE \|\hat T_{nm} - T_0\|_{L^2(\hat P_n)}^2 
\lesssim \bbE \Big[W_2^2(\hat P_n,P) + 
 W_2^2(\hat Q_m,Q)\Big].$$ 
In the smooth setting, we hypothesize that the distributions $P$ and $Q$ have densities which are smooth. We then need to provide bounds on the risk of $\hat P_n$ and $\hat Q_n$ as estimates of $P$ and $Q$ in the Wasserstein distance. Bounds of this type are well-studied in the literature for a variety of smooth density estimators~\citep{weed2019a,divol2022}. 
We then need to relate the $L^2(\hat P_n)$ 
risk to the $L^2(P)$ risk. The approach we follow 
is to show that our estimator $\hat P_n$
has density ratio $dP / d\hat P_n$ that is 
almost surely bounded over the supports of $P$ and $\hat P_n$, 
in which case the $L^2(P)$ risk is bounded
above by the $L^2(\hat P_n)$ risk up to constants.

 The following result appears in the work of~\citet{balakrishnan2025stability} who build on previous results of~\citet{manole2024plugin}.
\begin{theorem}[Informal]
\label{thm:smooth_densities}
For any $s > 0$, suppose that $p,q$ are $s$-smooth and strictly positive on a known regular domain $\Omega$. Then 
there exist proper and absolutely continuous  estimators 
$\hat P_n$ and $\hat Q_n$
such that the unique optimal
transport map $\hat T_n$ pushing
$\hat P_n$ forward onto $\hat Q_m$
satisfies
\begin{align} 
\label{eq:smooth_rate} 
\mathbb{E} \|\hat T_{nm} - T_0\|_{L^2(P)}^2 \leq 
C  \epsilon_{n\wedge m}  ,\quad \text{where } 
\epsilon_n:=\begin{cases} 
1/n, & d = 1, \\ 
\log n/n, & d = 2, \\
n^{-\frac{2(s+1)}{2s+d}}, & d \geq 3.
\end{cases}
\end{align} 
\end{theorem}

\noindent The convergence rate
appearing in Theorem~\ref{thm:smooth_densities}
matches known minimax lower bounds for estimating OT maps between H\"older-continuous densities~\citep{hutter2021}, up to a logarithmic factor when $d=2$.
One interesting aspect of these rates is that they are \emph{faster} than the usual rates for non-parametric density estimation in the $L^2$-norm which scale as $n^{-\frac{2s}{2s + d}}$. One intuition for this stems from the stability bound in Theorem~\ref{thm:stability} which upper bounds the convergence rates for the transport map estimate by the squared Wasserstein distance between the density estimate and truth. Under the assumption that the estimates and true distributions have densities which are bounded away from $0$ on their support, a result of \citet{peyre2018} shows that $$W_2(\widehat{P},P) \lesssim \|\widehat{p} - p\|_{\dot H^{-1}}:= \sup_{g: \|g\|_{\dot H^1} \leq 1} \big| \bbE_{\hat P}[g] - \bbE_{P} [g]\big|,$$
where $\|g\|_{\dot H^1} = \|\nabla g\|_{L^2}$
is the first-order (homogeneous) Sobolev
norm of a map $g$. 
This shows that the rate of convergence of $\widehat{P}$ to $P$ in the Wasserstein distance is upper bounded by the density estimation rate in an inverse Sobolev norm. The inverse Sobolev norm is a more forgiving metric than the usual $L^2$-norm studied in non-parametric density estimation -- it downweighs discrepancies in the higher frequency components of the difference $\widehat{p} - p$ -- and this leads to faster rates of convergence.

Despite its statistical optimality the smooth OT map can be difficult to compute in multivariate settings. One strategy to approximate the smooth OT map is to resample the density estimators $\widehat{P}_n$ and $\widehat{Q}_m$, and then to compute the nearest neighbor map between the resamples. This procedure will retain the attractive statistical properties of the plugin estimator but would require us to draw many more than $n \wedge m$ samples from the density estimates.
Quantitative convergence
rates have also been derived
for several other methods, 
based on
entropic  OT with vanishing regularization~\citep{pooladian2021,divol2024tight,eckstein2024convergence,mordant2024},
and sum-of-squares relaxations~\citep{vacher2021dimension,muzellec2021}. 
Although these methods 
are not known to be minimax optimal
in the same level of generality
as plugin or dual estimators, they typically enjoy 
more favorable computational properties.

\subsection{Inference for Smooth OT Maps}
\label{sec:clt_smooth}
We now consider the study of limit laws for the transport map under smoothness following closely the work of \citet{manole2023}. \citet{manole2023} develop a pointwise CLT for the OT map between $P$ and $Q$ under the assumption that the measures are smooth and are supported on the $d$-dimensional flat torus $\mathbb{T}^d$. Developing multivariate CLTs beyond this restrictive setting is an interesting direction for future work.

Suppose that we consider the one-sample problem, where the source measure $P$ is known, and only the target measure $Q$ is sampled. We obtain samples $Y_1,\ldots,Y_n \sim Q$, and construct the kernel density estimator:
\begin{align*}
    \widehat{q}(x) = \frac{1}{n} \sum_{i=1}^n K_h (Y_i - x),
\end{align*}
for an appropriate choice of the kernel $K_h$ and its bandwidth $h$. We denote by $\widehat{Q}$ the distribution with density $\widehat{q}$, and let $\widehat{T}$ denote the OT map from the distribution $P$ to $\widehat{Q}$. The main result of \citet{manole2023} is the following:
\begin{theorem}[Informal]
\label{thm:limit}
Let $d \geq 3$, and suppose that $p,q$ are $s$-smooth with $s > 2$ and strictly positive. 
Suppose that 
$$h = o(n^{-\frac 1 {2s+d}}),\quad 
\text{and,}\quad n^{\frac{1}{d+4}}h \to \infty.$$
Then, for any $x_0 \in \mathbb{T}^d$,
there exists
a positive definite
matrix $\Sigma(x_0)\in \bbR^{d\times d}$
such that
\begin{align*}
    \sqrt{nh^{d-2}} (\widehat{T}(x_0) - T_0(x_0)) \rightarrow N(0, \Sigma(x_0)).
\end{align*}
\end{theorem}
\noindent This result shows that the estimator $\widehat{T}(x_0)$ obeys a central limit theorem centered at its
population counterpart $T_0(x_0)$, when the bandwidth $h$ belongs to a specified range. The limiting variance can be consistently estimated from the samples, and one may even use the bootstrap for practical inference, at least in this idealized setting.

To provide some insight into the above result, we informally describe the strategy of \citet{manole2023}.
To begin, notice that an equivalent
characterization of the pushforward
condition ${T_0}_\# P = Q$
is for the densities
$p$ and $q$ of the measures $P$ and $Q$
to satisfy the {\it change-of-variable}
formula
$$p- \det(D T_0)q(T_0) =0 ,\quad \text{over } \bbT^d,$$
where $DT_0$ is the Jacobian of $T_0$. 
Rewriting this equation in terms of a potential
$\varphi_0$ satisfying $T_0 = \nabla\varphi_0$, 
we have 
\begin{equation} 
\label{eq:MA}
\Psi(\varphi_0) = 0, \quad \text{where } \Psi(\varphi):= p- \det(\nabla^2  \varphi)q(\nabla\varphi).
\end{equation}
Equation~\eqref{eq:MA} is a second-order
partial differential equation (PDE),
which falls into a family
of nonlinear equations known
as Monge-Amp\`ere equations~\citep{figalli2017},
and whose study has been at the heart
of many regularity properties
of OT maps in the smooth setting~\citep{dephilippis2014}.
The estimator $\hat T = \nabla\hat\varphi$  
also satisfies its own
Monge-Amp\`ere equation, 
which can be viewed as a peturbation of the one above:
 \begin{equation} 
\label{eq:MA_hat}
\hat \Psi(\hat \varphi) = 0,\quad \text{where } 
\hat\Psi(\varphi) := p- \det(\nabla^2  \varphi)\hat q(\nabla\varphi).
\end{equation}
This setting is reminiscent of $Z$-estimation, which defines an estimator as the solution of a system of ``estimating equations'' \citep{vandervaart1996}.
In this case, the estimator $\hat T$
is the gradient of the solution 
to an estimating equation, and its 
analysis
follows along similar conceptual lines -- by \emph{linearizing} the nonlinear estimating equation. It is straightforward to compute the Fr\'{e}chet derivatives $\Psi^{'}$ and $\widehat{\Psi}^{'}$ of the preceding operators. 
Ignoring remainder terms we conclude that,
\begin{align*}
    \widehat{\Psi}^{'}_{\varphi_0} (\widehat{\varphi} - \varphi_0) \approx \widehat{\Psi}(\widehat{\varphi}) - \widehat{\Psi}(\varphi_0), 
\end{align*}
where $\hat T = \nabla\hat\varphi$.
Continuing to reason informally, under appropriate conditions, we can replace the operator $\widehat{\Psi}^{'}_{\varphi_0}$ by $\Psi^{'}_{\varphi_0}$.
Furthermore, noting
that $\widehat{\Psi}(\widehat{\varphi}) - \widehat{\Psi}(\varphi_0) = \Psi(\varphi_0) - \widehat{\Psi}(\varphi_0) = \text{det}(\nabla^2\varphi_0) (\widehat{q}(\nabla \varphi_0) - q(\nabla \varphi_0))$, we obtain that:
\begin{align*}
\Psi^{'}_{\varphi_0} (\widehat{\varphi} - \varphi_0) \approx 
    \text{det}(\nabla^2\varphi_0) (\widehat{q}(\nabla \varphi_0) - q(\nabla \varphi_0))=:f.
\end{align*}
Let us abbreviate the operator
$\Psi^{'}_{\varphi_0}$
by $E$. Although we do not give the general
form of the operator
$E$ here, we mention that
$E$ is again a second-order partial differential
operator, but unlike $\Psi$, $E$ is {\it linear}.
Thus, the above display
tells us that $\hat\varphi - \varphi_0$ is
well-approximated by the solution to 
a linear, second-order PDE
of the form
\begin{equation}\label{eq:pde} 
Eu = f, \quad \text{over } \bbT^d.
\end{equation}
For example, 
when $P$ and $Q$ are uniform distributions
over the torus, $E$ turns out to be
the negative Laplacian operator $E = -\Delta = -\mathrm{trace}(\nabla^2(\cdot))$,
and the above equation
is known as 
the {\it Poisson equation}. 
The operator $E$ is
invertible,
thus one has the approximation
\begin{align*}
  \widehat{\varphi} -  \varphi_0 \approx E^{-1} f.
\end{align*}
Since $E$ is a linear operator, 
its inverse is also linear, 
and by definition of $f$, one 
can deduce that the entire expression
$E^{-1} f$ is a linear operator
in $\hat q-q$, which we abbreviate by $F(\hat q-q)$. Assuming we can take gradients
of the above display, we thus finally
arrive at the linearization
\begin{equation} 
\label{eq:linearization} 
\hat T - T_0 \approx \nabla F(\hat q-q).
\end{equation}
This suggests that in order to study (pointwise) limit laws for the OT map, we need to understand the (pointwise) fluctuations of the quantity $\nabla F(\widehat{q} - q)$, which is 
a linear smoother. 
This 
program is carried out rigorously by \citet{manole2023} who make all of these steps precise.
We also note
that the linearization~\eqref{eq:linearization}
has been generalized beyond
the torus to general Euclidean
domains in the work of~\cite{gonzalez2024linearization}. 
Statistical inference for related
PDE problems have also been
studied in the recent literature on 
infinite-dimensional 
inverse problems (e.g.~\cite{nickl2020convergence}). 

It is worth commenting on the relationship, in the smooth setting, between the stability bounds in Theorem~\ref{thm:stability} and the characterization used to establish the limit law in Theorem~\ref{thm:limit}. From the linearization~\eqref{eq:linearization}, one would expect that $\mathbb{E} \|\widehat{T} - T_0\|_{L^2(P)}^2 \approx \mathbb{E} \| \nabla F (\widehat{q} - q) \|_{L^2(P)}^2$. On the other hand, 
our discussion in the previous
section implied that $\mathbb{E} \|\widehat{T} - T_0\|_{L^2(P)}^2  \lesssim \mathbb{E} \|\widehat{q} - q\|^2_{\dot H^{-1}}$.
One can show that these two upper bounds
are on the same order.
Therefore, the pointwise characterization~\eqref{eq:linearization} 
implies results that are comparable to the stability bound in Theorem~\ref{thm:stability}. On the other hand, the stability bound only crudely (up to constants) upper bounds the global $L^2(P)$ error of the transport map, while the limit law builds on an asymptotically \emph{exact} characterization of the pointwise error. Obtaining this stronger, exact characterization is however more challenging and uses restrictive assumptions.

The idea of linearizing the Monge-Amp\`ere
equation has
also recently been used with great success 
in the probability theory literature,
for the problem of
deriving precise asymptotics
for the empirical in-sample
coupling $T_n$ described
in Section~\ref{sec:NN}. 
Unlike the plugin estimator $\hat T$
described here, the empirical map $T_n$
is highly nonsmooth, and the Monge-Amp\`ere
equation~\eqref{eq:MA_hat} does
 not make sense for this estimator. 
Nevertheless, 
starting with the work of~\cite{ambrosio2019b}, 
it was shown
in a series of papers
that the first-order
behavior of $T_n$ is governed by solutions to  equation~\eqref{eq:pde} in a 
weak sense,
leading to new asymptotic
results for the empirical 
map and the Wasserstein distance (e.g.~\cite{goldman2022,ambrosio2019,ambrosio2022,goldman2021quantitative,huesmann2024}). 
We expect that
these ideas could prove fruitful
to study OT map estimators in low-regularity settings.

\subsection{Beyond Smooth OT Maps}
\label{sec:beyond_smooth}
The results which we have surveyed in this section---and indeed much of
the existing literature on statistical inference for OT
maps between absolutely continuous distributions---has relied regularity conditions for $T_0=\nabla\varphi_0$,
typically in the form of smoothness and strong convexity conditions
on the potential $\varphi_0$. 
One can gauge the generality of these modeling assumptions  
in one of two ways. On the one hand, given {\it any} distribution $P$
and {\it any} strongly convex and smooth map $\varphi$, Brenier's theorem implies that
$\nabla \varphi$ is automatically the OT map from $P$ to $Q = \nabla\varphi_\# P$. 
This argument implies that $P$ could be arbitrarily irregular, and the conditions on $\varphi_0$ would hold
for appropriate choices of the target $Q$.
That is, regularity of the measures is not a prerequisite for the regularity
of the OT map. 
On the other hand, one can also state sufficient conditions on $P$
and $Q$ under which $T_0$ will be regular; 
one set of sufficient conditions are that
$P$ and $Q$ admit continuous densities supported on convex sets
with smooth boundary~\citep{caffarelli1992,caffarelli1996,dephilippis2014}.
We have also seen another
set of sufficient conditions in Section~\ref{sec:NN} for log-concave distributions~\citep{caffarelli2000,carlier2024optimal}.

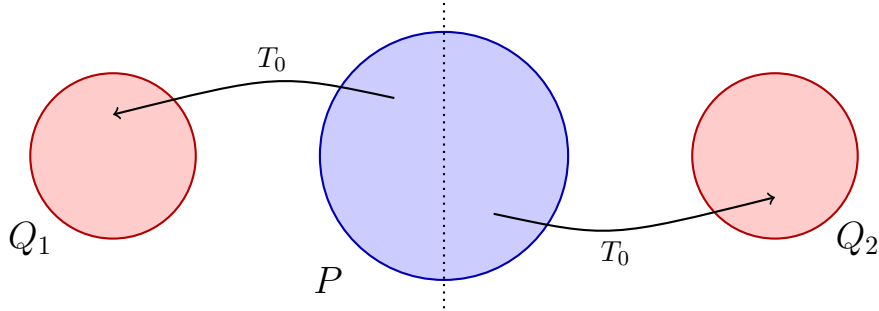
\begin{figure}[t] \centering 
\begin{tikzpicture}[scale=1.1]
\draw[fill=blue!20, draw=blue!70!black, thick] (0,0) circle (1.5);
\node at (-1.4,-1.5) {\Large $P$};

\draw[dotted, thick] (0,1.85) -- (0,-1.85);

\begin{scope}[shift={(-4,0)}]
  \draw[fill=red!20, draw=red!70!black, thick] (0,0) circle (1);
  \node at (-1,-1) {\Large $Q_1$};
\end{scope}

\begin{scope}[shift={(4,0)}]
  \draw[fill=red!20, draw=red!70!black, thick] (0,0) circle (1);
  \node at (1,-1) {\Large $Q_2$};
\end{scope}

 
\draw[->, thick] (-0.6,0.7) .. controls (-2,1) .. (-4,0.5)
      node[midway, above] {$T_0$};

\draw[->, thick] (0.6,-0.7) .. controls (2,-1) .. (4,-0.5)
      node[midway, below] {$T_0$};
\end{tikzpicture}
\caption{\label{fig:map}Example of a discontinuous OT map. The source
distribution $P$ is a uniform distribution over the blue ball, and the target distribution
$Q = \frac 1 2 (Q_1 + Q_2)$ is a uniform distribution over the union of the two red balls. 
The OT map $T_0$ from $P$ to $Q$ has a discontinuity along the dotted black line.} 
\end{figure}

Despite the generality of these conditions, 
there are many benign examples of measures $P$
and $Q$ for which the OT map is irregular. 
A simple way of generating such examples is to pick $P$ to have a convex
support, and $Q$ 
to have a support which is topologically disconnected (and hence nonconvex), 
as illustrated in Figure~\ref{fig:map}. 
OT maps between such distributions are always discontinuous. Although their set
of discontinuities has Lebesgue measure zero~\citep{figalli2010}, 
they still violate the global smoothness assumptions 
which are required by the results stated throughout
this section\footnote{Assumptions on the connectedness of the supports of $P$ and $Q$
arise in inferential problems for the Wasserstein distance~(e.g.~\cite{delbarrio2019a,staudt2025uniqueness}).}. 
For instance, one can derive simple counter-examples
to the stability bound in Theorem~\ref{thm:stability}
when $T_0$ is not globally smooth~\citep{merigot2020}.
Inferential results for OT maps 
remain elusive in this general regime, 
aside from some special cases discussed in Section~\ref{sec:variants} below.

Nevertheless, weaker variants of the stability bound in Theorem~\ref{thm:stability}
exist even without any smoothness assumptions on the OT maps, and such inequalities
can be used to obtain upper bounds on the minimax estimation risk under fewer assumptions. 
The study of such general quantitative stability bounds was initiated by the work of~\cite{berman2021}, 
and is now a highly active research area in the mathematical literature on optimal transport theory~\citep{merigot2020,delalande2023,letrouit2024,mischler2024,divol2024tight,carlier2024,kitagawa2025}. 
We refer to the forthcoming lecture notes of~\cite{letrouit2025} 
for a comprehensive survey.
As an example, the following one-sample quantitative stability bound is due to~\cite{delalande2023}. 
\begin{theorem}[Informal]
\label{thm:stability_general}
Let $P$ have a positive density over a connected domain with sufficiently regular boundary. 
Let $Q,\hat Q$ be arbitrary compactly-supported distributions. Then, the optimal transport
maps $T_0$ and $\hat T$ which respectively push $P$ forward onto $Q$ and $\hat Q$ satisfy
\begin{align*}
\|\hat T - T_0\|_{L^2(P)}^2 \lesssim W_2^{\frac 1 3}(\hat Q,Q).
\end{align*}
\end{theorem}
\noindent Theorem~\ref{thm:stability_general} can be viewed as a one-sample variant of Theorem~\ref{thm:stability}
in which no assumptions are placed on $T_0$, but instead some very mild regularity is placed
on the source measure\footnote{It can also be shown that OT maps are qualitatively stable 
under essentially no conditions on the source measure~\citep{segers2022}.} 
$P$. For example, Theorem~\ref{thm:stability_general} can be applied to the discontinuous map $T_0$ appearing
in Figure~\ref{fig:map}. The price to pay for this level of generality is the slower exponent $1/3$, as compared
to the exponent 2 in Theorem~\ref{thm:stability}. Although the sharpness of this exponent remains an open question, 
Theorem~\ref{thm:stability_general} already leads to upper bounds on the minimax estimation risk
of nonsmooth optimal transport maps in the one-sample setting. For example, if $T_m$ denotes
the empirical optimal transport map, pushing forward a sufficiently regular measure
$P$ onto an empirical measure $Q_m$, comprised of
iid samples $Y_1,\dots,Y_m$ from any measure $Q$, then Theorem~\ref{thm:stability_general} readily implies
$$\bbE  \| T_m-T_0\|_{L^2(P)}^2 \lesssim \bbE \big[W_2^{\frac 1 3}(Q_m,Q)\big] \lesssim m^{-\frac 1 {3d}},\quad \text{for all } d \geq 5.$$\
\indent The fact that OT maps are consistently estimable without any regularity
assumptions stands in  contrast to other familiar function estimation problems in nonparametric statistics---like
density estimation or nonparametric regression---where some regularity is typically needed for consistent estimation. 
The reason for this distinction is the fact that OT maps are {\it multivariate monotone} maps, by virtue of being 
gradients of convex functions. It is well-known that monotone function estimation problems---like CDF estimation
or isotonic regression---can be carried out without further structural assumptions, and this benefit of monotonicity
is shared by OT maps. Some works~\citep{delara2021} have leveraged this perspective
to build piecewise constant and monotone estimators of the OT map, which are similar in spirit to the empirical CDF.

\section{Variants and Special Cases}
\label{sec:variants}
Despite much interest, and much recent progress on the general multivariate OT problem, the statistical results require strong regularity conditions and the obtained results often suffer strongly from the curse of dimensionality. 
In this section we consider variants and special cases of OT which paint a more optimistic picture for statistical inference.

\subsection{The One-Dimensional Case}
\label{sec:oned}
In the one-dimensional case, when the source distribution has no atoms, the OT map from $P$ to $Q$ is $T_0(x) = F_Q^{-1}(F_P(x)),$
and this suggests the natural plugin estimator:
\begin{align}
\label{eqn:empirical}
    \widehat{T}(x) = \widehat{F}_Q^{-1}(\widehat{F}_P(x)),
\end{align}
where $\widehat{F}_P$ and $\widehat{F}_Q$ are the empirical CDFs, of the source and target samples. The study of the empirical CDF and quantile function is classical \citep{shorack1986empirical}, 
and their composition in equation~\eqref{eqn:empirical}
is sometimes known as the quantile-quantile (QQ) process,
for which strong approximation results are classical~\citep{aly1986strong}.
We will survey here a more recent treatment of the estimator $\hat T$, which is due to \citet{ponnoprat2024uniform}. 

If the distributions of $P$ and $Q$ are well-behaved we would expect that $\widehat{T}$ will inherit the convergence rate of the empirical CDF and quantile processes, and therefore converge to $T_0$ at parametric rates (both pointwise and in $L^2(P)$). 
Suppose that there is an interval $[a,b]$ such that, (i) $F_Q$ is continuously differentiable on the interval $[F_Q^{-1}(F_P(a)),F_Q^{-1}(F_P(b))]$, (ii) the density of $Q$, $F_Q'$ is non-zero on the interval $[F_Q^{-1}(F_P(a)),F_Q^{-1}(F_P(b))]$. Then for any $x \in [a,b]$
\citet{ponnoprat2024uniform} derive the following Bahadur representation:
\begin{align*}
    \sqrt{n+m}(\widehat{T}(x) - T_0(x)) = \frac{\sqrt{n+m}}{F_Q'(T_0(x))} \left[ \frac{1}{n} \sum_{i=1}^n \mathbb{I}(X_i \leq x) + \frac{1}{m} \sum_{j=1}^m \mathbb{I}(Y_j \leq T_0(x)) \right] + o_{P \otimes Q}(1).
\end{align*} 
Although we are not aware of a reference, it is possible to derive a non-asymptotic analogue of the result of \citet{ponnoprat2024uniform} to show that under appropriate conditions on the source and target distributions, the empirical estimate~\eqref{eqn:empirical} converges at the parametric $1/n \wedge 1/m$ rate in the (squared) $L^2(P)$ metric. 
The Bahadur representation also allows one to derive a pointwise Normal limit and naturally enables bootstrap-based (pointwise) statistical inference. Concretely, fixing an $x \in [a,b]$ one may estimate the $1-\alpha/2$ quantile $\kappa(1-\alpha/2)$
for the quantity $\sqrt{n+m}(\widehat{T}^*(x) - \widehat{T}(x))$, where $\widehat{T}^*$ is the empirical transport map computed on bootstrap samples. Then, as shown by \citet{ponnoprat2024uniform} the interval 
$C_{\alpha} = \left[\widehat{T}(x) \pm \frac{\kappa(1-\alpha/2)}{\sqrt{n + m}}\right]$ is an asymptotically valid $1-\alpha$ confidence interval for the OT map.

Going beyond pointwise inference \citet{ponnoprat2024uniform} show that kernel smoothing enables uniform inference for the transport map. They show that under appropriate assumptions,
\begin{align*}
    \sqrt{n+m} \sup_{x \in [a,b]} \frac{\widehat{T}(x) - T_0(x)}{\widehat{\sigma}_x} \overset{d}{\rightarrow} \sup_{x \in [a,b]} |Z(x)|,
\end{align*}
where $Z$ is a mean-zero, Gaussian process, for an appropriate plugin standard deviation estimate $\hat\sigma_x$. Here $\widehat{T}(x)$ is taken to be a plugin estimate of the form in~\eqref{eqn:empirical}, compute with kernel smoothed empirical measures.
\citet{ponnoprat2024uniform} further show that the bootstrap confidence band is consistent at level $1-\alpha.$ 

The main takeaway from these results is that in the one-dimensional setting under relatively transparent assumptions, the OT map can be estimated at parametric rates, and practical inference is possible via the bootstrap.

\subsection{The Gaussian Case}
When $P$ and $Q$ are multivariate Normal, the OT map 
once again has a closed form expression in~\eqref{eqn:gaussianot} and one can construct
the natural plugin estimate:
\begin{align*}
    \widehat{T}(x) = \widehat{\Sigma}_1^{-1/2} (\widehat{\Sigma}_1^{1/2} \widehat{\Sigma}_2 \widehat{\Sigma}_1^{1/2} )^{1/2} \widehat{\Sigma}_1^{-1/2} (x - \widehat{\mu}_1) + \widehat{\mu}_2,
\end{align*}
where each of the estimates $\widehat{\mu}_1, \widehat{\mu}_2$, $\widehat{\Sigma}_1, \widehat{\Sigma}_2$ are constructed in the usual way.
The analysis of this estimator builds on the large body of recent work in non-asymptotic high-dimensional statistics. 
For a matrix $\Sigma$, we can define its effective rank as $r(\Sigma) = \text{tr}(\Sigma)/\lambda_{\max}(\Sigma)$.
We assume that the population covariances are well-conditioned, i.e. that:
\begin{align*}
    c < \lambda_{\min}(\Sigma_1), \lambda_{\min}(\Sigma_2) \leq \lambda_{\max}(\Sigma_1), \lambda_{\max}(\Sigma_2) < C, 
\end{align*}
for finite constants $c, C > 0$. Further, assume that $n_1 \gtrsim d$ (which ensures that $\widehat{\Sigma}_1$ is invertible with probability 1), and that $n_2 \gtrsim r(\Sigma_2)$. Then \citet{flamary2019concentration} show that with probability at lesat $1 - \exp(-t) - 1/n_1$
\begin{align*}
\|\widehat{T} - T_0\|_{L^2(P)} \lesssim \sqrt{r(\Sigma_1)} \times \left( \sqrt{\frac{r(\Sigma_1)}{n}} \vee \sqrt{\frac{r(\Sigma_2)}{m}} \vee \sqrt{\frac{t}{n \wedge m}} \vee \frac{t}{n \wedge m} \right).
\end{align*}
The rates for transport map estimation in this setting depend on the dimension primarily through the effective rank of the population covariance matrices, highlighting that these transport map estimates can be useful even in high-dimensional settings where the data dimension $d$ is comparable to the sample-size $n$. Under fixed-dimensional asymptotics, one can also derive a central limit theorem for $\hat T$ via the delta method \citep{vandervaart1998}. 

More broadly, imposing parametric structure either on the distributions $P$ and $Q$, on the transport map $T_0$, or on both, are an effective way to attempt to circumvent the curse of dimensionality and enable reasonable statistical inference even from moderate sample-sizes. 

\subsection{Discrete and Semi-Discrete Optimal Transport}

Our emphasis thus far has been 
on the study of OT maps between absolutely
continuous distributions. 
One of the downsides
of this setting is that, 
absent any strong structural
assumptions---such as the Gaussianity
assumption from the preceding
section---the minimax rate
of estimating OT maps degrades exponentially
with the ambient dimension $d$
(cf. Section~\ref{sec:NN}). 

One of the striking findings of the 
recent statistical optimal transport
literature has been the so-called
{\it lower complexity adaptation} (LCA) principle~\citep{forrow2019statistical,hundrieser2024lca,stromme2024minimum,groppe2024lower,delbarrio2024central}, which
asserts that the minimax rate
of estimating various OT objects
improves when  {\it at least}
one of the two measures $P$ and $Q$
is intrinsically low-dimensional.
That is, the fundamental hardness
of these estimation problems is driven
by the {\it least} complex of the two measures. 
This phenomenon was first crystallized by
the work of~\cite{hundrieser2024lca} in the context of estimating
Wasserstein distances.
Roughly speaking, their work
proves  that
if $P$ and $Q$ are measures
on $\bbR^d$, but one of them is
supported on a set which is intrinsically
$s$-dimensional, then 
the Wasserstein distance $W_2(P,Q)$
can be (adaptively) estimated at a rate which 
depends only on $s$ and not on $d$. 
This remarkable property significantly
expands the scope of applications
where the Wasserstein
distance can be tractably estimated. 

It has been conjectured
that the LCA principle
extends to the problem
of estimating optimal transport maps and couplings, 
and although a general theory of this
type has yet to be developed, 
the special case of {\it semi-discrete OT}
has been well-studied. 
One version of this setting 
is to assume that the source
distribution $P$ is 
absolutely
continuous, while the target distribution $Q$
is discrete with a finite number of atoms. In this setting, 
$Q$ morally has zero intrinsic dimension, 
and if the LCA principle
were to apply, one would
expect to be able to estimate
the OT map $T_0$
at a dimension-free rate,
which indeed turns out to be the case. 
As an example,
the work of~\cite{sadhu2024stability}
considers the following setting:
Assume $P$ is known, and that
one has i.i.d. samples
$Y_1,\dots,Y_m \sim Q$, where $Q$ is discrete and supported on a finite number of atoms.
A natural estimator of the OT map $T_0$
is the empirical map introduced
in Section~\ref{sec:NN}, namely
the OT map which pushes
forward $P$ onto the empirical measures
$Q_m = (1/m) \sum_{i=1}^m \delta_{Y_i}$. 
Under mild
assumptions on $P$,~\cite{sadhu2024stability}
show that the process 
$$\sqrt m\|T_m - T_0\|_{L^p(P)}^p$$
admits a non-degenerate limiting
distribution, for any fixed
$p \in [1,\infty)$.
This result 
proves that the squared 
$L^2(P)$ risk of the empirical OT map
decays at the dimension-free rate $n^{-1/2}$.
In contrast, we saw
in Section~\ref{sec:NN}
that the typical convergence rate of this same
estimator is
$n^{-2/d}$ when~$Q$
is absolutely continuous,
which depends exponentially on the dimension $d$. 
Furthermore, the result
of~\cite{sadhu2024approximation}
 enables the construction
of $L^p$-confidence bands
for the semi-discrete optimal
transport map under very mild assumptions, 
which again stands in contrast
to the absolutely
continuous setting, where
we saw that considerably stronger
conditions are needed to perform inference. 

Inferential
results can also be derived 
when both $P$ and~$Q$ are discrete---a setting
which was recently considered in the
work of~\cite{sommerfeld2018}. 
In this case, there generally
does not exist an OT map from $P$ to $Q$, 
and existing results
aim to provide inferential
methods for an optimal coupling
$\pi$ from $P$ to $Q$. 
One important subtlety here
is the fact that the coupling
$\pi$ need not be unique; nevertheless \cite{klatt2022limit}
derived central limit theorems
for empirical optimal couplings which are
centered at a sequence of optimal
couplings between $P$ and $Q$. 
These results were followed up by the work
of~\cite{liu2023asymptotic}, who derived 
a tractable procedure for building
confidence sets based on the asymptotics
of~\cite{klatt2022limit}, 
and~\cite{liu2025beyond} who derive limit
laws for optimal couplings.

\subsection{Entropic Optimal Transport}
The entropic OT problem can be studied as an approximation to the classical OT problem, as the regularization parameter $\varepsilon \rightarrow 0$. This perspective can lead to various computational
benefits \citet{pooladian2021,altschuler2017near,cuturi2013}, but does not 
lead to statistical benefits in the estimation of classical OT quantities. On the other hand, when the regularization parameter $\varepsilon$ is held fixed, and the distributions $P$ and $Q$ are compactly supported the resulting entropic potentials, maps and distances can be estimated at fast parametric rates.

Focusing on the entropic transport map~\eqref{eqn:entropicmap}, the work of \citet{rigollet2025sample}, building on earlier work by \citet{mena2019statistical,genevay2019sample} shows that for compactly supported $P$ and $Q$:
\begin{align*}
    \mathbb{E} \|\widehat{T}_{\varepsilon} - T_{\varepsilon}\|_{L^2(P)}^2 \lesssim \frac{1}{n},
\end{align*}
where the implicit constant depends exponentially on $\varepsilon$. We remark that this situation mirrors the situation in much of non-parametric statistics. For instance, the \emph{smoothed} density (say, the density convolved with a Gaussian of bandwidth $h$) can be estimated at fast parametric rates under mild assumptions, provided the bandwidth is held fixed. Whether the entropic map (or the smoothed density in non-parametrics) is a reasonable target depends on the application domain and on the downstream use of the estimator.

Recent work has also shown that the benefits of entropic OT also extend to statistical inference \citep{gonzalez-sanz2022weak,goldfeld2024limit,gonzalez2023weak,harchaoui2020asymptotics,gunsilius2021matching}. For instance, the work of \citet{goldfeld2024limit} shows that, when the distributions are compactly supported, for any fixed $s \in \mathbb{N}$:
\begin{align}
\label{eqn:process}
\sqrt{n} (\widehat{T}_{\varepsilon} - T_{\varepsilon}) \overset{d}{\rightarrow} - \nabla G,
\end{align}
where $G$ is a mean zero Gaussian random variable which takes value in a space of $s$-smooth functions. This result should be contrasted with the unregularized, classical OT setting, where only pointwise limits have been studied and it is known that there is no scaling of the process $\widehat{T} - T_0$ which has a non-degenerate limit (where $\widehat{T}$ is the kernel density plugin estimate with a reasonable bandwidth choice) \citep{manole2023}. Building on the result in~\eqref{eqn:process}, \citet{goldfeld2024limit} also show that the bootstrap can be used to construct valid \emph{confidence bands} for the entropic OT map (as in Section~\ref{sec:oned}).
Taken together these results highlight some of the many statistical benefits of the entropic OT framework.

\subsection{Divergence-Regularized Optimal Transport}
One possible drawback of entropic optimal transport is the fact that it leads to {\it dense} optimal couplings,
which are always supported on the product
space $\supp(P) \times \supp(Q)$. 
This behavior stands in contrast to conventional optimal transport, for which the optimal coupling is supported on the degenerate
set $\{(x,T_0(x)): x \in \supp(P)\}$, at least when 
there exists an optimal transport map $T_0$ from $P$ to $Q$. 
This overspreading effect of entropic OT is one of the prices to pay for its 
favorable smoothness properties.

A series of recent works has shown that one can achieve a middle ground between these two behaviors, using
a variant of entropic optimal transport known as {\it divergence-regularized} optimal transport~\citep{essid2018quadratically,blondel2018smooth,muzellec2017tsallis,dessein2018regularized}.
The divergence-regularized OT problem consists of replacing the Kullback-Leibler
penalization in equation~\eqref{eqn:eot}  by a penalty involving other $f$-divergences,
such as Tsallis divergences, which include the $\chi^2$-divergence. 
It has been observed empirically that such regularization procedures typically lead to sparse
couplings, and these findings have recently been complemented by theoretical results in some special cases~\citep{garriz2024infinitesimal,wiesel2024sparsity,gonzalez2024sparsity}. 
For example, in the one-dimensional case,~\cite{gonzalez2024sparsity} derive quantitative 
bounds on a notion of width of the support of the $\chi^2$-regularized
optimal transport coupling, showing that it decays polynomially in the regularization parameter~$\varepsilon$. 

Although the  divergence-regularized OT problem does not typically enjoy the same smoothness
properties as entropic OT~\citep{bayraktar2025stability}, it was shown by \cite{gonzalez2025sparse} that, much like entropic OT, 
certain Tsallis-regularized couplings can be estimated at the parametric rate, independent of the ambient dimension.
\cite{gonzalez2025sparse}~also show that these regularized couplings obey a central limit theorem
when integrated against smooth test functions, much like in our discussion of entropic OT~\eqref{eq:test_function_clt}. 
These results suggest that the statistical properties of certain divergence-regularized OT 
couplings mirror that of entropic OT, despite its sparsity-inducing properties.

\section{Discussion and Future Work}
\label{sec:discussion}
Despite exciting recent progress on the estimation of and inference for transport maps, several important open problems remain. For the estimation of the OT map, a better understanding would be achieved by studying estimation without strong curvature conditions, and under alternative structural assumptions on the measures. One concrete direction is to investigate the tightness of the statistical rates that are derived from Theorem~\ref{thm:stability_general}.

This review has focused primarily on the OT map for the quadratic cost, where the theory is most mature. Broadening our statistical understanding to encompass other cost functions is a natural next step. The analysis of plug-in estimators typically imposes assumptions on the measures $P$ and $Q$ and on the OT map $T_0$, while dual estimators are analyzed under conditions on $T_0$ alone. Developing techniques to close this gap and analyze plug-in estimators using only properties of $T_0$ would be interesting. In contrast to the nearest neighbor map, the smooth plug-in OT map estimators described in Section~\ref{sec:smooth} face practical limitations: na\"{i}ve computation can require resampling the density estimates prohibitively many times, followed by computing the nearest neighbor map on the resamples. More practical estimators that retain the statistical advantages of the smooth plug-in estimators would be a valuable contribution.

A complementary line of work would involve jointly analyzing the computational and statistical properties of OT map estimation methods. While many numerical procedures exist for computing the OT map when the distributions are known, studying their behavior in the statistical setting—where access to the distributions is stochastic—remains largely unexplored. Joint characterizations of iteration complexity and statistical accuracy could enable a practitioner to meaningfully tradeoff between computational efficiency and statistical performance.

For inference on the multivariate OT map, a key open problem is to establish pointwise limit laws under more realistic assumptions. Central limit theorems for the $\ell_\infty$ error in high-dimensional settings are also largely absent. Progress in this area would require a deeper, global understanding of the limiting behavior of OT map estimators.

Entropic OT has long been known to have good computational properties~\citep{cuturi2013}. Recent works have highlighted that these benefits also extend to estimation and inference.
One important drawback is that entropic OT is not always natural, and more work needs to be done to suggest data-driven choices for the regularization parameter $\varepsilon$. As we noted earlier, a similarly optimistic statistical picture would arise in essentially all of non-parametric statistics if we targeted fixed-bandwidth smoothed estimands. Arguably in non-parametric statistics the theory and practice of tuning-parameter selection is better developed, and this is another area where we expect to see more work in the future.

Beyond optimal transport and entropic optimal transport there are many natural ways to construct transport maps. These include the Knothe-Rosenblatt rearrangement~\citep{rosenblatt1952remarks,knothe1957contributions,baptista2025conditional}, diffusion models~\citep{ho2020denoising,song2021score}, rectified flow~\citep{liu2023flow}, sliced OT~\citep{bonneel2015sliced,manole2022b} and Gaussian smoothed optimal transport~\citep{goldfeld2020convergence}.
The statistical theory of these alternative ways to construct transport maps is in its infancy and we hope to see significant progress studying these methods, their statistical properties, and the tradeoffs in using them for various applications.

\section*{Acknowledgments}
The authors are grateful to Jonathan Niles-Weed for many discussions which formed the basis of this review. 
TM also thanks Axel Munk, Shayan Hundrieser, Cyril Letrouit, Alberto Gonz\'alez-Sanz, and Marcel Nutz, for 
discussions which helped shape this survey. This work was supported in part by the NSF grant DMS-2310632.
TM is supported by a Norbert Wiener postdoctoral fellowship.

\bibliographystyle{abbrvnat}
\bibliography{map_review}

\appendix

\end{document}